\documentclass[Journal,9pt]{IEEEtran} % versione technical note
\usepackage{amssymb,bm}
\usepackage{makeidx}         % allows index generation
\usepackage{amsmath}
\usepackage{multicol}        % used for the two-column index
\usepackage[bottom]{footmisc}% places footnotes at page bottom
\usepackage{graphicx}
\usepackage{color}
\usepackage{amsfonts}
\usepackage[latin1]{inputenc}
\usepackage{epstopdf}
\usepackage[usenames,dvipsnames]{xcolor}
\usepackage[ruled,vlined]{algorithm2e}

\newcommand{\eq}{\begin{equation}}
\newcommand{\eeq}{\end{equation}}
\newcommand{\eqn}{\begin{eqnarray}}
\newcommand{\eeqn}{\end{eqnarray}}

\newcommand{\bsea}{\begin{subeqnarray}}
\newcommand{\esea}{\end{subeqnarray}}
\newcommand{\nn}{\nonumber}

\newcommand{\al}[1]{\begin{align}  #1 \end{align}}\newcommand{\xib}{\bm \xi}

\def\bmat{\left[ \begin{array}}
\def\emat{\end{array} \right]}
\newcommand{\Range}{\mathop{\rm Range}}

\newcommand{\tr}{\mathop{\rm tr}}  %traccia

\newcommand{\Cc}{ \mathcal{C}}

\newcommand{\Fc}{ \mathcal{F}}

\newcommand{\Jc}{ \mathcal{J}}

\newcommand{\Lc}{ \mathcal{L}}

\newcommand{\Pc}{ \mathcal{P}}
\newcommand{\Qc}{ \mathcal{Q}}

\newcommand{\Sc}{ \mathcal{S}}

\newcommand{\Cs}{ \mathbb{C}}

\newcommand{\Es}{ \mathbb{E}}

\newcommand{\Ls}{ \mathbb{L}}

\newcommand{\Rs}{ \mathbb{R}}

\newcommand{\Ts}{ \mathbb{T}}

\newcommand{\Zs}{ \mathbb{Z}}

\newcommand{\mz}{\color{black}}

\def\qed{\hfill \vrule height 7pt width 7pt depth 0pt \smallskip}

\newcounter{pippo}

\newtheorem{remark}{Remark}[section]
\newtheorem{teor}{Theorem}[section]
\newtheorem{corr}{Corollary}[section]
\newtheorem{propo}{Proposition}[section]
\newtheorem{lemm}{Lemma}[section]
\newtheorem{exam}{Example}
\newtheorem{probl}[pippo]{Problem}
\newtheorem{defn}{Definition}[section]

\newcommand{\teo}{\begin{teor}}
\newcommand{\eteo}{\end{teor}}
\newcommand{\cor}{\begin{corr}}
\newcommand{\ecor}{\end{corr}}
\newcommand{\prop}{\begin{propo}}
\newcommand{\eprop}{\end{propo}}
\newcommand{\lem}{\begin{lemm}}
\newcommand{\elem}{\end{lemm}}
\newcommand{\ex}{\begin{exam}}
\newcommand{\eex}{\end{exam}}
\newcommand{\pb}{\begin{probl}}
\newcommand{\epb}{\end{probl}}
\newcommand{\df}{\begin{defn}}
\newcommand{\edf}{\end{defn}}
\newcommand{\aprop}{\begin{apropo}}
\newcommand{\eaprop}{\end{apropo}}
\newcommand{\alem}{\begin{alemm}}
\newcommand{\ealem}{\end{alemm}}
\newcommand{\rem}{\begin{remark}}
\newcommand{\erem}{\end{remark}}

\newcommand{\xb}  {\mathbf{x}}
\newcommand{\yb}  {\mathbf{y}}
\newcommand{\eb}  {\mathbf{e}}

\begin{document}

\title{\mz Optimal Transport between Gaussian Stationary Processes}

\author{Mattia~Zorzi \thanks{M. Zorzi is with the
Dipartimento di Ingegneria dell'Informazione, Universit\`a degli studi di
Padova, via Gradenigo 6/B, 35131 Padova, Italy
({\tt\small zorzimat@dei.unipd.it}). This work was partially supported by the SID project ``A Multidimensional and Multivariate Moment Problem Theory for Target Parameter Estimation in Automotive Radars'' (ZORZ\_SID19\_01) funded by the Department of Information Engineering of the  University of Padova. }}

\markboth{DRAFT}{Shell \MakeLowercase{\textit{et al.}}: Bare Demo of IEEEtran.cls for Journals}

\maketitle

\begin{abstract} We consider the optimal transport problem between multivariate Gaussian stationary stochastic processes. The transportation effort is the variance of the filtered discrepancy process. {\mz The main contribution of this technical note is to show} that the corresponding solution leads to a weighted Hellinger distance between multivariate power spectral densities. Then, we propose a spectral estimation approach in the case of indirect measurements which is based on this distance. 
 \end{abstract}


\section{Introduction}
The optimal mass transport problem recently gained attention in the control filed. Its modern formulation, due to {\em Kantorovich} \cite{kantorovich1942translocation}, consists in minimizing the effort of transporting one nonnegative measure to another nonnegative given an associated cost of moving mass from a point to another one. 
The latter has been used to derive new distances between covariance matrices and spectral densities \cite{georgiou2008metrics,elvander2018interpolation}. Matrix-valued transport problems have been proposed in \cite{ning2014matrix,chen2017matricial} which lead to transportation distances between multivariate power spectral densities. The key idea in these works is to set up the optimal transport problem in terms of the joint power spectral density. It is worth noting that the latter  gives a complete joint description of two zero-mean multivariate Gaussian stationary stochastic {\mz processes}.

In the present paper, instead, we face the optimal transport problem in a different way. More precisely, such a problem is formulated {\mz directly} in terms of the joint probability density of two Gaussian stationary processes. The transportation cost is the variance of the filtered discrepancy process, i.e. the difference between the two processes. {\mz The main contribution of this technical note is to} show that the corresponding solution leads to a weighted {\em Hellinger} distance between multivariate power spectral densities. The latter has been introduced in \cite{Hellinger_Ferrante_Pavon,ramponi2009globally} for the case without weight function. It is worth noticing that the aforementioned  correspondence can be regarded as the dynamic extension of to the results showed in \cite{NING_2013,knott1984optimal} for zero-mean Gaussian random vectors.

Then, we show that the weighted transportation distance can be used to perform spectral estimation in the case we have indirect measurements. More precisely, we have a network of $m$ sensors which measure indirectly $m$ variables of interest for which we want to estimate the spectral density. The matrix-valued weight function of our distance is given by the transfer matrix characteristic of the sensor network. We consider the THREE-like spectral estimation paradigm for which a large body of literature for the scalar \cite{A_NEW_APPROACH_BYRNES_2000,georgiou2003kullback,BYRNES_GUSEV_LINDQUIST_RATIONAL_COV_EXT,karlsson2013uncertainty,ALPHA}, multivariate \cite{FERRANTE_TIME_AND_SPECTRAL_2012,BETAPRED,BETA,8706531,zhu2018existence} and multidimensional case \cite{georgiou2006relative,SIAM_IM_COMPR,ringh2018multidimensional,zhu2019m} has been developed. Here, the estimated spectral density is the   closest one to a prior spectral density and matching the output covariance of a bank of filters (fed with the measured data). The prior represents the priori information that we have about the spectrum. The bank of filter is designed by the user in such a way the estimator exhibit high-resolution properties in prescribed frequency bands, \cite{A_NEW_APPROACH_BYRNES_2000}. 
Finally, we compare the proposed estimator with the one obtained by using the {\em Itakura-Saito distance}. The latter has been in introduced in \cite{FERRANTE_TIME_AND_SPECTRAL_2012}.

The outline of the technical note is as follows. In Section 
\ref{sect_Transp} we introduce the optimal transport problem between  Gaussian processes and we derive the corresponding weighted Hellinger distance. Section \ref{sect_spectral_est} is devoted to spectral estimation with indirect measurements. In Section \ref{sect_results} we provide  a simulation study to test the proposed estimator. Finally, in Section \ref{sect_concl} we draw the conclusions.

The following notation will be adopted throughout the paper. $\Qc_n$ denotes the space of symmetric matrices of dimension $n\times n $.   Given an Hermitian matrix $X$, $X>0$ ($X\geq 0$) means that $X$ is positive (semi-)definite.  Given $X\geq 0$, $X^{1/2}$ denotes a square root of $X$, i.e. $X=X^{1/2}X^{1/2}$. Given a matrix $X\in\Cs^{n\times n}$: $[X]_{il}$ denotes its entry in position $(i,l)$; we define the norms $\|X\|: =\sqrt{\tr(X^*X)}$ and $\|X\|_W:=\sqrt{\tr(X^* W X)}$ with $W=W^*>0$. We define the unit circle as $\Ts=\{e^{j\vartheta} \hbox{ s.t. } \vartheta\in[0,2\pi)\}$. $\Ls_\infty^{m\times m}(\Ts)$ denotes the 
function space of matrix functions defined on $\Ts$ and taking values in $\Cs^{n\times n}$.
Given $\Phi,\Psi\in\Ls_\infty^{m\times m}(\Ts)$: the shorthand notation $\int_\Ts \Phi$ denotes the integration of $\Phi$ on $\Ts$ with respect to the normalized Lebesgue measure; $\Phi>0$ means that $\Phi(e^{j\vartheta})>0$  for any $e^{j\vartheta}\in\Ts$ almost everywhere,  $\Phi=\Psi$ means that $\Phi$ and $\Psi$ coincide almost everywhere. $\Sc_m^+$ denotes the set of spectral densities bounded and coercive of dimension $m\times m$. Given a sequence $h=\{ h_t,\; t\in\Zs\}$, $(h)_t$ denotes $h_t$. {\mz Let $g$ and $h$ denote two sequences, then $g\star h$ denotes the discrete convolution operation.}

\section{Transportation distance between Gaussian processes}\label{sect_Transp}
We consider two $m$-dimensional jointly Gaussian stationary stochastic processes $\xb=\{\xb_t,\; t\in\Zs\}$ and $\yb=\{\yb_t,\; t\in\Zs\}$ whose mean is equal to zero. The latter are completely described by the finite dimensional probability density functions $p_\xb(x_t,x_s; t,s)$ and $p_\yb(y_t,y_s;t,s)$ with $t,s\in\Zs$. Let $p_{\xb,\yb}(x_t,x_s,y_u,y_v;t,s,u,v)$, with $t,s,u,v\in\Zs$, be the finite dimensional joint probability density of $\xb$ and $\yb$. Then, we consider the following optimal transport problem:
\al{\label{pb_trasp}d(p_\xb,p_\yb)^2= \underset{p_{\xb,\yb}\in \Pc}{\inf} \{\Es[\| (h\star (\xb- \yb))_{t}\|^2]  \hbox{ s.t. (\ref{cond_marg1})-(\ref{cond_marg2}) hold}  \} }
where $\Pc$ denotes the set of Gaussian joint probability densities $p_{\xb,\yb}$; $h=\{ h_t,\; t\in\Zs\}$, with $h_t\in\Rs^{m\times m}$ and $h_t=0$ for any $t<0$, represents the impulse of a BIBO system; 
\al{\label{cond_marg1}&\int_{\Rs^m}\int_{\Rs^m} p_{\xb,\yb}(x_t, x_s,y_u,y_v;t,s,u,v)\mathrm d y_u \mathrm d y_v= p_{\xb}(x_t, x_s;t,s), \\
&\label{cond_marg2} \int_{\Rs^m}\int_{\Rs^m} p_{\xb,\yb}(x_t,x_s,y_u , y_v;t,s,u,v)\mathrm d x_t \mathrm d x_s= p_{\yb}(y_u,y_v;u,v),\ } for any $t,s,u,v\in\mathbb Z$. It is worth noting that (\ref{pb_trasp}) represents the optimal transport problem between Gaussian processes $\xb$ and $\yb$. The transportation cost is the variance of the filtered process $h\star (\xb-\yb)$ and $\xb-\yb$ can be regarded as the discrepancy process. 

{\em Example}. In the case that $h$ is such that $h_0=I$ and $h_t=0$  for any $t\neq 0$, the transportation cost is the variance of the discrepancy process $\xb-\yb$. Then, we have
\al{d(p_\xb,p_\yb)=\left[\underset{p_{\xb,\yb}\in \Pc}{\inf} \{\Es[\| \xb_t- \yb_{t}\|^2]  \hbox{ s.t. (\ref{cond_marg1})-(\ref{cond_marg2}) hold}  \}\right]^{1/2}.\nn} The latter is the dynamic counterpart of the 2-Wasserstein distance in \cite{NING_2013} between Gaussian random vectors.

{\mz  Let $\Phi_\xb$ and $\Phi_\yb$ denote the power spectral density of $\xb$ and $\yb$, respectively. The next result shows that the solution of the optimal transport problem in (\ref{pb_trasp}) can be expressed in terms of $\Phi_\xb$ and $\Phi_\yb$. \prop Consider the optimization problem in (\ref{pb_trasp}). It holds that:
\al{\label{def_uitle}d_\Omega(p_\xb,p_\yb)^2&=\tr\int_\Ts \Omega\Phi_\xb +\Omega\Phi_\yb -2(\Phi_\yb^{1/2} \Omega \Phi_\xb \Omega \Phi_\yb^{1/2}) ^{1/2}.}\eprop 
\IEEEproof} Consider a joint process $[\xb^T\,\yb^T]^T$ whose probability density belongs to $\Pc$ and satisfying (\ref{cond_marg1})-(\ref{cond_marg2}).  Since the latter is Gaussian, an equivalent description for such a process is given by its  power spectral density \al{\Phi=\left[\begin{array}{cc}\Phi_\xb & \Phi_{\xb \yb}\\ \Phi_{\yb\xb} & \Phi_\yb\end{array}\right].\nn} {\mz Notice that}
\al{\Es[(h\star \xb)_t (h\star \xb)_t^T]&=\int_\Ts H\Phi_\xb H^*\nn\\
\Es[(h\star \yb)_t (h\star \yb)_t^T]&=\int_\Ts H\Phi_\yb  H ^*\nn\\
\Es[(h\star \xb)_t (h\star \yb)_t^T]&=\int_\Ts H \Phi_{\xb\yb} H ^*\nn}
where $H$ denotes the discrete-time Fourier transform of $h$:\al{\label{def_H}H(e^{j\vartheta})=\sum_{t\geq 0} h_t e^{-j\vartheta t}.}   We can rewrite the objective function in (\ref{pb_trasp}) in terms of $\Phi$:
\al{  \Es & [ \| (h  \star(\xb-\yb))_t\|^2] =\tr \Es[ (h\star(\xb-\yb))_t(h\star(\xb-\yb))_t^T]\nn\\  & =\tr\Es[(h\star \xb))_t(h\star \xb)_t^T+
(h\star \yb)_t(h\star \yb)_t^T\nn\\
&\hspace{1cm}-(h\star \xb)_t(h\star \yb)_t^T -(h\star \yb)_t(h\star \xb)_t^T  ]\nn\\
&= \tr \int_\Ts H (\Phi_\xb +\Phi_\yb   -\Phi_{\xb \yb} -\Phi_{\yb \xb} )  H^* \nn\\
&= \tr\int_\Ts \Omega (\Phi_\xb +\Phi_\yb - \Phi_{\xb \yb} - \Phi_{\yb \xb} )       \nn }
where $\Omega(e^{j\vartheta})  =  H(e^{j\vartheta})^*   H(e^{j\vartheta})$. In what follows we assume that $\Omega,\Phi_\xb,\Phi_\yb\in\Sc_m^+$. The latter condition is not so restrictive. Indeed, it means those spectral densities admit a minimum phase spectral factor. Therefore, we can rewrite (\ref{pb_trasp}) in terms of $\Phi $:
\al{ d_\Omega(p_\xb,p_\yb)^2= &\underset{\Phi_{\xb\yb}  }{\inf}   \tr\int_\Ts \Omega (\Phi_\xb +\Phi_\yb - \Phi_{\xb \yb} -\Phi_{\yb \xb} )     \nn\\ &    \hbox{ s.t. }  \left[\begin{array}{cc}\Phi_\xb & \Phi_{\xb \yb}\\ \Phi_{\yb\xb} & \Phi_\yb\end{array}\right]   \geq  0 \nn}
 where the infimization is only with respect to $\Phi_{\xb\yb}$ because $\Phi_\xb$ and $\Phi_\yb$ have been already fixed by constraints  (\ref{cond_marg1})-(\ref{cond_marg2}). Notice that we added the subscript $\Omega$ to $d$ in order to stress such a dependence. Since $\Phi_\yb\in\Sc_m^+$, we have $\Phi_\yb >0$. Accordingly, condition $\Phi\geq 0 $ is equivalent to $\Phi_\xb - \Phi_{\xb \yb} \Phi_\yb^{-1} \Phi_{\yb\xb}\geq 0$ and hence
\al{\label{pb_trasp2} d_\Omega(p_\xb,p_\yb)^2= &\underset{\Phi_{\xb\yb}  }{\inf}   \tr\int_\Ts \Omega (\Phi_\xb+\Phi_\yb- \Phi_{\xb \yb}  - \Phi_{\yb \xb} )       \nn\\ &    \hbox{ s.t. }  \Phi_\xb - \Phi_{\xb \yb} \Phi_\yb^{-1} \Phi_{\yb\xb}     \geq  0.}
 Let \al{\Lambda:=\Phi_\xb - \Phi_{\xb \yb} \Phi_\yb^{-1} \Phi_{\yb\xb},\nn } then
\al{(\Phi_\xb-\Lambda)^{1/2}\Upsilon =\Phi_{\xb \yb} \Phi_\yb^{-1/2}\nn } where $\Upsilon$ is an all-pass function, i.e. $\Upsilon (e^{j\vartheta})\Upsilon (e^{j\vartheta})^*=I$ for any $e^{j\vartheta}\in\Ts$. Thus, \al{\Phi_{\xb\yb}= (\Phi_\xb-\Lambda)^{1/2 } \Upsilon \Phi_\yb^{1/2}.\nn} Substituting $\Phi_{\xb\yb}$ with $\Lambda$ and $\Upsilon$ in (\ref{pb_trasp2}) we obtain 
\al{ \label{pb_Lamb}d_\Omega(p_\xb,p_\yb)^2-\int_{\Ts}\Omega(\Phi_\xb+\Phi_\yb)&=\underset{\Lambda }{\,\inf\,} \Fc(\Lambda)\,\nn\\
&\hspace{0.5cm} \hbox{ s.t. }  \Lambda\geq0 }
where 
\al{\label{pb_L_fix}\Fc(\Lambda)=&\underset{\Upsilon}{\,\inf\,}  \tr  \int_\Ts (  \Psi\Upsilon  +   \Upsilon^* \Psi ^*)    \nn\\ &    \hspace{0.1cm}\hbox{ s.t. }   \Upsilon \Upsilon^*=I}
 and \al{\label{def_Psi}\Psi= -\Phi_\yb^{1/2} \Omega (\Phi_\xb-\Lambda)^{1/2}.}
 We solve problem (\ref{pb_L_fix}) by means of duality theory. The Lagrangian is 
 \al{L(\Upsilon, \Delta )&= \tr  \int_\Ts  (  \Psi \Upsilon  +   \Upsilon^* \Psi^*)    +\tr\int_\Ts (\Upsilon  \Upsilon^*-I) \Delta }
where $\Delta\in \Ls_\infty^{m\times m}(\Ts)$ and such that $\Delta=\Delta^*$ is the Lagrange multiplier. Moreover, the first and the second variation of $L(\cdot, \Delta)$ along $\delta\Upsilon \in \Ls_\infty^{m\times m}(\Ts)$ are
\al{\delta L(\Upsilon, \Delta  ; \delta\Upsilon)&= \tr\int_\Ts \Psi\delta \Upsilon +\delta\Upsilon^* \Psi^*+ \delta \Upsilon \Upsilon^* \Delta+  \Upsilon\delta   \Upsilon^* \Delta    \nn\\
\delta^2 L(\Upsilon, \Delta  ; \delta\Upsilon)&=2 \tr\int_\Ts  \delta   \Upsilon^* \Delta \delta   \Upsilon\nn.}
Accordingly, $L$ is lower bounded if and only if $\Delta\geq 0$.
Moreover, if $\Delta>0$ then $L$ is strictly convex with respect to $\Upsilon$ and the unique point of minimum is given by the stationarity condition $\delta L(\Upsilon,\Delta;\delta\Upsilon)=0$ for any $\delta\Upsilon \in \Ls_\infty^{m\times m}(\Ts)$ which implies the optimal form \al{\Upsilon^\circ = -\Delta^{-1}\Psi^*.\nn}
Accordingly, the dual functional is
\al{\Jc(\Delta)= L(\Upsilon^\circ,\Delta)=-\tr  \int_\Ts   \Psi \Delta^{-1}\Psi ^*+\Delta \nn} and the dual problem is $\max_{\Delta>0} \, J(\Delta)$. Notice that the first and the second variation along $\delta \Delta\in\Ls_\infty^m(\Ts)$ are 
\al{\delta \Jc(\Delta ; \delta \Delta)&= \tr\int_\Ts \Psi \Delta^{-1} \delta\Delta \Delta^{-1}\Psi ^*+\delta\Delta\nn\\ 
\delta^2 \Jc(\Delta ; \delta \Delta)&= -2\tr\int_\Ts \Psi \Delta^{-1} \delta\Delta \Delta^{-1} \delta\Delta \Delta^{-1}\Psi ^*. \nn }
Under the assumption that  $\Delta>0$ and $\Psi$ is full rank almost everywhere, we have that $\delta^2 \Jc(\Delta ; \delta \Delta)<0$ for any $\delta\Delta\in \Ls_\infty^{m\times m}(\Ts)$ and $\delta\Delta\neq 0$. Accordingly, $\Jc$ is strictly concave, and the unique point of minimum is given by setting equal to zero its first variation in any direction which gives
\al{\Delta^2=\Psi^*\Psi.\nn} Therefore the point of maximum is 
\al{\Delta^\circ =(\Psi^*\Psi)^{1/2}.\nn} Notice that $\Delta^\circ >0$ (i.e. our assumption on $\Delta$ is satisfied) provided that $\Psi$ is full rank almost everywhere. In this case, the unique solution to (\ref{pb_L_fix}) is \al{\Upsilon^\circ=-(\Delta^\circ)^{-1}\Psi^*=-(\Psi^*\Psi)^{-1/2}\Psi^*.\nn} Moreover, it is not difficult to see that  \al{\label{optUPsPsi}\Psi \Upsilon^\circ=-(\Psi \Psi^*)^{1/2}=-(\Phi_\yb^{1/2} \Omega (\Phi_\xb-\Lambda) \Omega \Phi_\yb^{1/2}) ^{1/2} } where we have exploited (\ref{def_Psi}). Substituting (\ref{optUPsPsi}) in (\ref{pb_L_fix}) we obtain:
\al{\Fc(\Lambda)&= 2\tr\int_{\Ts}\Psi\Upsilon^\circ = -2\tr\int_\Ts(\Phi_\yb^{1/2} \Omega (\Phi_\xb-\Lambda) \Omega \Phi_\yb^{1/2}) ^{1/2} .\nn} Accordingly, the optimal solution to (\ref{pb_Lamb}) is $\Lambda^\circ=0$ and thus: 
{\mz \al{d_\Omega(p_\xb,p_\yb)^2&=\tr\int_\Ts \Omega\Phi_\xb +\Omega\Phi_\yb -2(\Phi_\yb^{1/2} \Omega \Phi_\xb \Omega \Phi_\yb^{1/2}) ^{1/2}\nn } which gives (\ref{def_uitle}).} Notice that $\Psi=-\Phi_\yb^{1/2}\Omega\Phi_\xb^{1/2}$ is full rank almost everywhere for $\Lambda^\circ=0$, i.e. our assumption on $\Psi$ is satisfied. {\mz \qed\\}

\prop For any square spectral factor $W_\xb$ of $\Phi_\xb$, we have 
\al{\label{second_def}d_\Omega(p_\xb,p_\yb)^2=\min    & \{\int_{\Ts}  \| W_\xb- W_\yb\|_{\Omega}^2 \hbox{ s.t. }  \nn \\  &  W_\yb\in\Ls_\infty^{m\times m}(\Ts),  W_\yb W_\yb^*=\Phi_\yb\}.} Moreover, $d_\Omega(p_\xb,p_\yb)$ is a bona fide distance function. \eprop
\IEEEproof Once $W_{\xb}$ is fixed, any square spectral factor $ W_{\xb}^\prime$ of $\Phi_\xb$ can be written as $ W_\xb^\prime=W_\xb \Upsilon$ where $\Upsilon\in\Ls_\infty^{m\times m}(\Ts)$ is an all-pass function. Therefore, we have 
\al{\tr\int_\Ts \|W_\xb^\prime-W_\yb\|^2_\Omega=\tr\int_\Ts (W_\xb^\prime-W_\yb)^*\Omega(W_\xb^\prime-W_\yb)\nn\\=\tr\int_\Ts (W_\xb-W_\yb\Upsilon^*)^*\Omega(W_\xb-W_\yb\Upsilon^* )\nn} 
and $W_\yb\Upsilon^*$ is a spectral factor a $\Phi_\yb$. Therefore, the right hand side of (\ref{second_def}) does not depend on the particular choice of the square spectral factor for $\Phi_\xb$. To prove the equality in (\ref{second_def}), observe that the right hand side can be rewritten as 
\al{\min \left\{\int_\Ts\| W_\xb-\Phi_\yb^{1/2} \Upsilon\|_\Omega^2, \; \Upsilon\in \Ls_\infty^{m\times m}  (\Ts), \; \Upsilon\Upsilon^*=I\right\}.\nn} The corresponding Lagrangian is 
\al{L(\Upsilon,\Lambda)=\int_\Ts\| W_\xb-\Phi_\yb^{1/2} \Upsilon\|_\Omega^2+\tr(\Upsilon\Upsilon^*-I)\Lambda\nn}
where $\Lambda=\Lambda^*\in\Ls_\infty^{m\times m}(\Ts)$ is the Lagrange multiplier. Under the assumption that $\Lambda>0$, the minimum of $L$ with respect to $\Upsilon$ is $\Upsilon(\Lambda)=\Lambda^{-1}\Phi^{1/2}_\yb\Omega W_\xb$. It is not difficult to see that the Lagrange multiplier satisfying the constraint in the primal problem is $\Lambda^\circ= \Phi_\yb^{1/2}\Omega \Phi_\xb\Omega \Phi_\yb^{1/2}>0$. Therefore, the minimum point of the right hand side in (\ref{second_def}) is 
\al{W_\yb^\circ=\Phi^{1/2}_\yb(\Phi_\yb^{1/2}\Omega \Phi_\xb \Omega\Phi_\yb^{1/2})^{-1/2}\Phi_\yb^{1/2}\Omega W_\xb\nn}
and by direct substitution we obtain (\ref{def_uitle}).   
The last part of the claim can be easily checked.  
\qed

In the special case that $\Omega=I$, we obtain the Hellinger distance for multivariate power spectral densities \cite{Hellinger_Ferrante_Pavon}:\al{d(p_\xb,p_\yb):=\left[ \inf \{\int_{\Ts}\| \right. & W_\xb  -W_\yb\|^2 \hbox{ s.t. }    \nn\\ & \left. W_\yb\in\Ls_\infty^{m \times m}(\Ts), W_\yb W_\yb^*=\Phi_\yb\}\right]^{1/2}.\nn} Accordingly, $d_\Omega(p_\xb,p_\yb)$ is the weighted Hellinger distance with weight the matricial function $\Omega\in\Sc_m^+$. 
Finally, notice that $d_\Omega(p_\xb,p_\yb)$ is completely characterized by $\Phi_\xb$ and $\Phi_\yb$. Accordingly, in what follows we will use the notation $d_\Omega(\Phi_\xb,\Phi_\yb)$ to denote the weighted transportation distance between the Gaussian processes having power spectral density $\Phi_\xb$ and $\Phi_\yb$, respectively. 

\section{Spectral Estimation with indirect measurements}\label{sect_spectral_est}
\begin{figure}
\centering
\includegraphics[width=0.7\columnwidth]{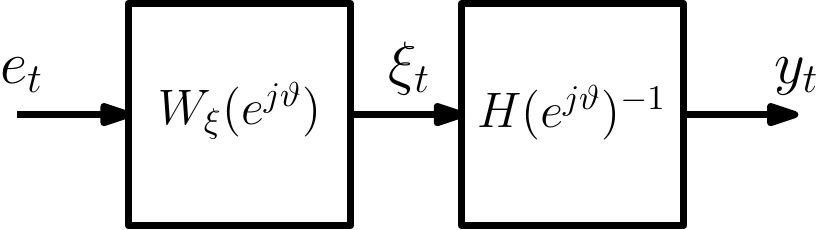}\caption{Generative model for processes $y$ and $\xi$.}\label{fig_process}
\end{figure}

{\mz In this section we show how the weighted Hellinger distance can be used to estimate stochastic processes from indirect measurements. As we will see, the appealing property is that the corresponding estimator accounts for the transfer matrix characteristic of the sensor network.} Consider two $m$-dimensional zero mean Gaussian stationary stochastic processes $\yb$ and $\xib$ generated by a normalized white Gaussian noise process $\eb$ as depicted in Figure \ref{fig_process}.  We assume that $H$, defined in (\ref{def_H}), is a rational transfer matrix, causally invertible and known. The transfer matrix $W_{\xib}$ is such that $\Phi_{\xib}=W_{\xib} W^*_{\xib}\in\Sc_m^+$ and it is not known. It is worth noting that the spectral density of $\yb$ is $\Phi_{\yb}=H^{-1} \Phi_{\xib} H^{-*}\in\Sc_m^+$ while the one of $\xib$ is $\Phi_{\xib}$. Next, we face the problem to find an estimate of $\Phi_{\xib}$ from a finite length realization $y^N=\{ y_1\ldots y_N\}$ of $\yb$. In plain words, we want to perform a multivariate spectral estimation task from indirect measurements: we have a network of $m$ sensors modeled by $\yb$ and measuring indirectly $m$ variables of interest modeled  by $\xib$. Moreover, $H^{-1}$ is the transfer matrix characteristic of the sensor network: $[H^{-1}]_{il}$ denotes the transfer function between the $l$-th variable of interest and the $i$-th sensor.

We address such a problem using a THREE-like framework. Suppose that a power spectral density (called prior) $\Psi_{\xib}\in\Sc_m^+$, which gives some a priori information on $\Phi_{\xib}$, is available. If no prior is available we choose the uninformative prior $\Psi_{\xib}=I$, i.e. the spectral density corresponding to white noise. Then, we fix a bank of filters \al{G(z)=(zI-A)^{-1}B\nn}
where $A\in\Rs^{n\times n}$ is a stability matrix, $B\in\Rs^{n\times m}$ is full column rank with $n>m$, and $(A,B)$ is a reachable pair. Let $\Sigma$ denote the steady state covariance of the output of the bank of filters 
\al{x_{t+1}=Ax_{t}+By_t.\nn} We compute an estimate $\hat \Sigma$, based on the data $y^N$, of $\Sigma$. Then, an estimate of $\Phi_{\xib}$ is given by the spectrum approximation problem:
 \al{\label{Pb_PSD_orig}\hat \Phi_{\xib,T}=& \underset{\Phi_{\xib}\in\Sc_m^+}{\mathrm{argmax\,}} d(\Phi_{\xib},\Psi_{\xib})\nn\\ &
\hbox{ s.t. } \Phi=H^{-1}\Phi_{\xib}H^{-*}, \;  \int_{\Ts} G\Phi G^*=\hat \Sigma} where ``T'' stand for transportation distance.    
The latter problem can be rewritten only in terms of $\Phi$:
 \al{\label{Pb_PSD}\hat \Phi=& \underset{\Phi\in\Sc_m^+}{\mathrm{argmax\,}} d_\Omega(\Phi,\Psi)\nn\\ &
\hbox{ s.t. } \int_{\Ts} G\Phi G^*=\hat \Sigma}
where $\hat \Phi_{\xib,T}=H\hat \Phi H^*$, $\Psi_{\xib}=H \Psi H^{*}$ and $\Omega=H^*H$. In the special case that $G(z)=[\,z^{-l}I_m \, \ldots \,z^{-1}I_m\,]^T$, we have that $\Sigma$ is a $n\times n$, with $n=lm$, block Toeplitz matrix whose first block row is $[\,\Sigma_0\, \Sigma_1 \,\ldots \,\Sigma_{n-1}\, ]$ and $\Sigma_k=\mathbb E[y_ty_{t+k}^T]$. Therefore, with this particular choice of $G(z)$, we fix the first covariance lags of $\hat \Phi$, i.e. $\hat \Phi$ is the solution to a covariance extension problem. It is worth noting $\hat \Sigma$ must guarantee the feasibility of Problem (\ref{Pb_PSD}) that is there exits at least one $\tilde \Phi\in\Sc_m^+$ such that $\int_{\Ts} G\tilde \Phi G^*=\hat \Sigma$. In \cite{georgiou2002structure}, it has been shown that $\hat \Sigma$ guarantees the feasibility of the problem if and only if $\hat \Sigma>0$ and $\hat \Sigma\in \mathrm{Range\,}\Gamma$ where $\Gamma$ is the linear operator defined as \al{\Gamma\,:\; &\Cc_m(\Ts) \rightarrow \Qc_n\nn\\
& \Phi\mapsto \int_\Ts G\Phi G^*\nn} where $\Cc_m(\Ts)$ denotes the family of $\Cs^{m\times m}$ and hermitian-valued continuous functions on the unit circle $\Ts$.  Optimization approaches have proposed to compute such an estimate from data, see \cite{zorzi2012estimation,5953479}.

 Let $W_\Psi$ and $W$ be a square spectral factor of $\Psi$ and $\Phi$, respectively, then Problem (\ref{Pb_PSD}) is equivalent to 
 \al{\label{Pb_PSD2}\hat W=& \underset{W\in\Ls_\infty^{m\times m}(\Ts)}{\mathrm{argmax\,}} \tr\int_{\Ts}(W- W_\Psi)(W-W_\Psi)^*\Omega\nn\\ &
\hbox{ s.t. } \int_{\Ts} GWW^*G^*=\hat \Sigma.}   
 The Lagrangian of Problem (\ref{Pb_PSD2}) is
\al{L(W,\Lambda)=&\tr\int_{\Ts} (W-W_\Psi)(W-W_\Psi)^*\Omega\nn\\ &+\tr\left(\Lambda\int_{\Ts} GWW^*G^*\right)-\tr(\Lambda\hat \Sigma)\nn\\
=&\tr\int_{\Ts} (\Omega+G^*\Lambda G)WW^*-\Omega W_\Psi W^*-\Omega WW_\Psi^*\nn\\ &+\Omega W_\Psi W_\Psi^*-\tr(\Lambda\hat \Sigma)\nn}
where $\Lambda\in \Qc_n$ is the Lagrange multiplier. $L(\cdot,\Lambda)$ is bounded below only if $\Omega+G^*\Lambda G>0$. In this case the infimum is attained and given by setting its first variation equal to zero
\al{\delta L(W,\Lambda;\delta W)=&\int_\Ts(\Omega+G^*\Lambda G)\delta WW^*+(\Omega+G^*\Lambda G)W\delta W^*\nn\\ &-\Omega W_\Psi \delta W^*-\Omega \delta WW_\Psi^* =0\nn}
for any $\delta W\in \mathbb L_\infty^{m\times m}(\Ts)$. The latter leads to the optimality condition
\al{(\Omega+G^*\Lambda G)W-\Omega W_\Psi=0\nn}
and thus the optimal form of $W$ is
\al{\label{opt_W}W^\circ=(\Omega+G^*\Lambda G)^{-1}\Omega W_\Psi.} Therefore, the dual problem consists in finding a point of maximum $\Lambda^\circ$ for 
\al{L(W^\circ,\Lambda)=-\left[\tr\int_{\Ts} \Omega\Psi \Omega (\Omega+G^*\Lambda G)^{-1}-\Omega\Psi\right]-\tr(\hat \Sigma\Lambda)\nn} over the set \al{\Lc_+=\{\,\Lambda\in\Qc_n \hbox{ s.t. } \Omega+G^*\Lambda G>0\}.\nn}
Thus, such a problem is equivalent to minimize 
\al{J(\Lambda)= \tr\int_{\Ts} \Omega\Psi \Omega (\Omega+G^*\Lambda G)^{-1}+\tr(\hat \Sigma\Lambda)\nn}
with $\Lambda\in\mathcal L_+$. As it has been shown in \cite{Hellinger_Ferrante_Pavon}, $G^*\Lambda G =0$ for any $\Lambda\in[ \Range\,\Gamma]^\bot$. Thus, $J(\Lambda+\Lambda_\bot)=J(\Lambda)$ for any $\Lambda_\bot\in[ \Range\,\Gamma]^\bot$. Accordingly, we can restrict the search of $\Lambda^\circ$ over the set $\Lc_{\Gamma,+}=\Lc_+ \cap \Range\, \Gamma$. Accordingly, we obtain the equivalent problem
\al{\label{pb_dual}\Lambda^\circ =\underset{\Lambda \in \Lc_{\Gamma,+}}{\mathrm{argmin\,}}J(\Lambda).} 

\teo \label{teo_esistenza}Problem (\ref{pb_dual}) admits a unique solution.\eteo
\IEEEproof Following arguments similar to the ones in \cite[Lemma 7.5]{Hellinger_Ferrante_Pavon}, it is possible to prove that $J$ is strictly convex on $\Lc_{\Gamma,+}$. Accordingly, if Problem (\ref{pb_dual}) admits solution, then the latter is unique. 

Next, we show the existence of such solution. First, 
notice that the set $\Lc_{\Gamma,+}$ is nonempty, indeed $0\in\Lambda_{\Gamma,+}$, open and unbounded. Moreover, $J$ is a continuous function on $\Lc_{\Gamma,+}$. The idea is to show that Problem (\ref{pb_dual}) is equivalent to minimize $J$ on a compact set, say $ \Lc^\star_{\Gamma,+}$. Then, by the {\em Weierstrass} Theorem we conclude that the minimum exists. 

We proceed to characterize $\Lc^\star_{\Gamma,+}$. First, we show that $J$ is bounded below on $\Lc_{\Gamma,+}$:
\al{J(\Lambda)&=\tr \int_{\Ts}\Psi^{1/2} \Omega (\Omega+G^*\Lambda G)^{-1}\Omega \Psi^{1/2}+\tr(\hat\Sigma\Lambda)\nn\\ 
&\geq \tr(\hat\Sigma\Lambda)=\tr\int_{\Ts} G\tilde \Phi G^* \Lambda=\tr\int_{\Ts} \tilde \Phi^{1/2}G^* \Lambda G\tilde \Phi^{1/2}\nn\\ &=
\tr\int_{\Ts} \tilde \Phi^{1/2} (\Omega+ G^* \Lambda G)\tilde \Phi^{1/2}-\tr\int_{\Ts}\tilde \Phi^{1/2}\Omega \tilde \Phi^{1/2} \nn\\
& \geq -\int_{\Ts}\tilde \Phi^{1/2}\Omega\tilde \Phi^{1/2}>-\infty\nn} where we exploited the following facts: $\Omega(e^{j\vartheta})+G(e^{j\vartheta})^*\Lambda G(e^{j\vartheta})$ is positive definite for any $e^{j\vartheta}\in\Ts$; $\Omega$ and $\tilde \Phi$ are bounded and coercive. Following arguments similar to the ones in \cite[Theorem 7.7]{Hellinger_Ferrante_Pavon}, it is possible to prove that 
\al{\underset{\Lambda\rightarrow \partial \Lc_{\Gamma,+}}{\lim}\, J(\Lambda)=\infty \nn\\
\underset{\|\Lambda\|\rightarrow \infty}{\lim}\, J(\Lambda)=\infty\nn} where $\partial \Lc_{\Gamma,+}$ denotes the boundary of $\Lc_{\Gamma,+}$ that is the set of $\Lambda\in\Lc_{\Gamma,+}$ such that $\Omega(e^{j\bar \vartheta})+G(e^{j\bar \vartheta})^*\Lambda G(e^{j\bar \vartheta})\geq 0$ and singular for at least one $e^{j\bar \vartheta}\in\Ts$.  Accordingly, we can restrict the search of the minimum point of $J$ on the  bounded and closed (and thus compact)  set $\Lc_{\Gamma,+}^\star=\{\,\Lambda\in\Lc_{\Gamma,+} \hbox{ s.t. } \Omega+G^*\Lambda G\geq \alpha I, \; \|\Lambda\|\leq \beta \}$ for some $\alpha,\beta>0$. \qed

\begin{corr} Problem (\ref{Pb_PSD}) admits the unique solution
\al{\label{opt_PHI}\hat \Phi_T=(\Omega+G^*\Lambda^\circ G)^{-1}\Omega\Psi\Omega(\Omega+G^*\Lambda^\circ G)^{-1}} where $\Lambda^\circ$ is the unique solution to (\ref{pb_dual}). \end{corr}
\IEEEproof Recall that the optimal solution $\Lambda^\circ$ of (\ref{pb_dual}) is a stationary point for $J$. The latter condition implies that  $W^\circ$ in  (\ref{opt_W}) with $\Lambda=\Lambda^\circ$ is such that the spectral density $W^\circ (W^\circ)^*$ satisfies the moment constraints in (\ref{Pb_PSD2}). Accordingly, the duality gap is equal to zero. Hence, Problem (\ref{Pb_PSD}) admits the unique solution $\hat\Phi_T=W^\circ (W^\circ)^*$ with $\Lambda=\Lambda^\circ$.\qed

Since $\Omega=H^*H$ and $H$ is causally invertible, then 
(\ref{opt_PHI}) can be rewritten as 
\al{\hat \Phi_T=H^{-1}&(I+H^{-*}G^*\Lambda^\circ GH^{-1})^{-1}\Psi_{\xib}\nn\\ &\times (I+H^{-*}G^*\Lambda^\circ GH^{-1})^{-1}H^{-*}\nn} 
and thus 
\al{\hat \Phi_{\xib,T}=(I+H^{-*}G^*\Lambda^\circ GH^{-1})^{-1}\Psi_{\xib}(I+H^{-*}G^*\Lambda^\circ GH^{-1})^{-1}.\nn}

An alternative to Problem (\ref{Pb_PSD_orig}) is to consider 
\al{\label{Pb_PSD_origIS}\hat \Phi_{\xib,IS}=& \underset{\Phi_{\xib}\in\Sc_m^+}{\mathrm{argmax\,}} d_{IS}(\Phi_{\xib}\|\Psi_{\xib})\nn\\ &
\hbox{ s.t. } \Phi=H^{-1}\Phi_{\xib}H^{-*}, \;  \int_{\Ts} G\Phi G^*=\hat \Sigma}   
where $d_{IS}$ is the Itakura-Saito distance:
\al{d_{IS}(\Psi_{\xib}\|\Phi_{\xib})=  \tr\int_{\Ts}(\log |\Psi_{\xib}| -\log |\Phi_{\xib}|+ \Phi_{\xib} \Psi_{\xib}^{-1}-I_m ).\nn} On the other hand, it is not difficult to see that
\al{d_{IS}(\Psi_{\xib}\|\Phi_{\xib})=d_{IS}(H^{-1}\Psi_{\xib}H^{-*}\|H^{-1}\Phi_{\xib}H^{-*})=d_{IS}(\Psi\|\Phi).\nn}
Accordingly, Problem (\ref{Pb_PSD_origIS}) is equivalent to solve the THREE-like estimator using direct measurements proposed in \cite{FERRANTE_TIME_AND_SPECTRAL_2012}:
\al{\label{Pb_PSD_origIS2}\hat \Phi_{IS}=& \underset{\Phi\in\Sc_m^+}{\mathrm{argmax\,}} d_{IS}(\Phi\|\Psi)\nn\\ &
\hbox{ s.t. }   \int_{\Ts} G\Phi G^*=\hat \Sigma}   
 where $\hat \Phi_{\xib,IS}=H\hat \Phi_{IS}H^*$.
 The optimal solution of (\ref{Pb_PSD_origIS2}) is
\al{\hat\Phi_{IS}&=(\Psi^{-1}+G^*\Lambda G)^{-1}
\nn\\
& = H^{-1}(\Psi_{\xib}^{-1}+H^{-*}G^*\Lambda GH^{-1})^{-1} H^{-*}\nn} where $\Lambda$ is given by solving the corresponding dual problem, see \cite{FERRANTE_TIME_AND_SPECTRAL_2012}.
Hence, we have
\al{\hat \Phi_{\xib,IS}=(\Psi_{\xib}^{-1}+H^{-*}G^*\Lambda GH^{-1})^{-1}.} \begin{figure}[ht]
\centering
\includegraphics[width=\columnwidth]{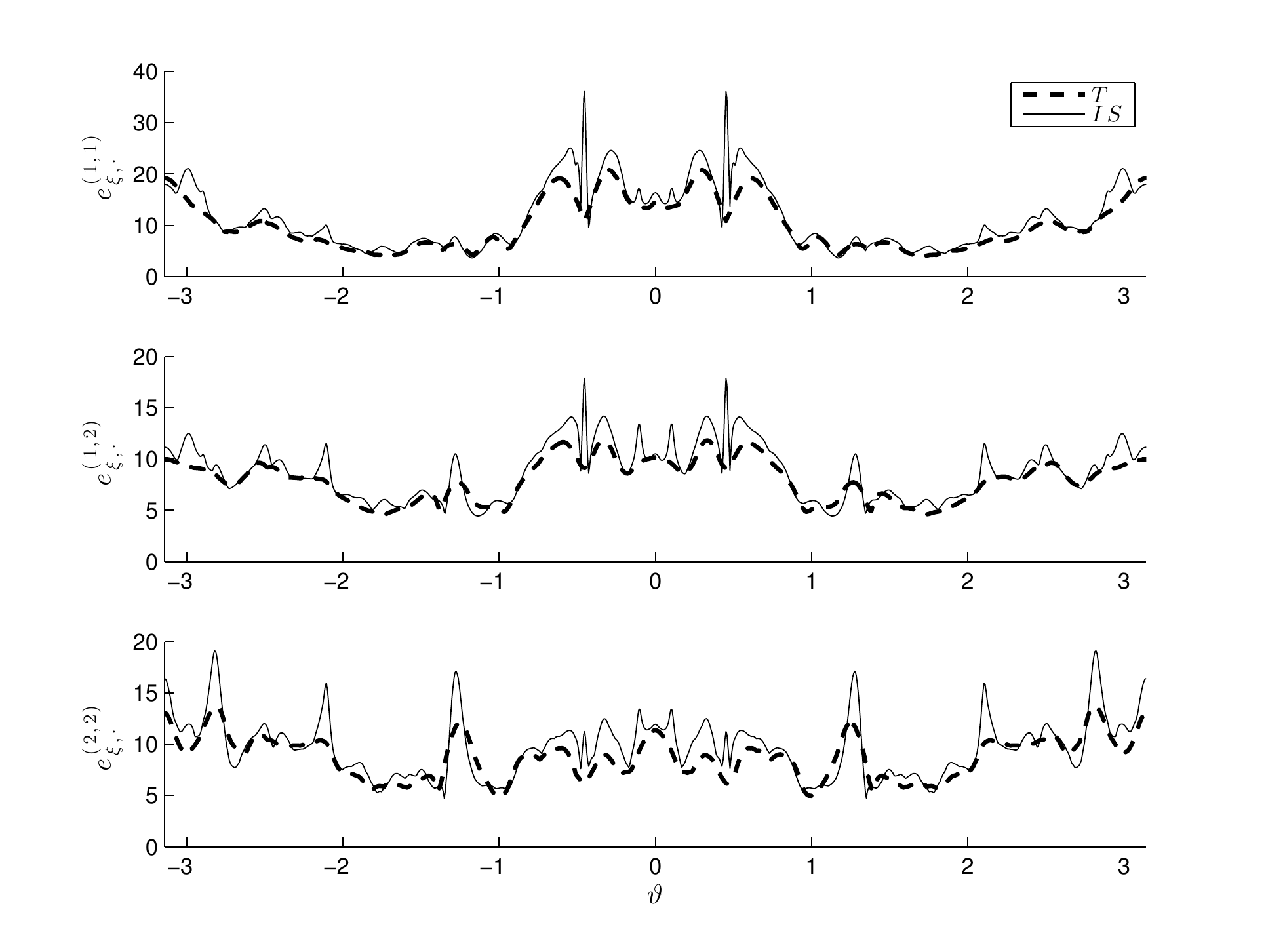}
\caption{Average errors in the estimation of $\Phi_{\xib}$ using the estimator with the transportation distance (T) and the Itakura-Saito distance (IS).}\label{fig_ris}
\end{figure} The fact that the objective function in (\ref{Pb_PSD_origIS2}) is the same for the cases of indirect/direct measurements reveals that the Itakura-Saito distance does not adapt the optimization scheme for a possible improvement in the case of indirect measurements. Drawing inspiration from (\ref{Pb_PSD}), we could consider the weighted Itakura-Saito distance, see \cite{DUAL}, with weight function $\Omega$:
\al{\label{Pb_PSD_IS}\hat \Phi=& \underset{\Phi\in\Sc_m^+}{\mathrm{argmax\,}} d_{IS,\Omega}(\Phi,\Psi)\nn\\ &
\hbox{ s.t. } \int_{\Ts} G\Phi G^*=\hat \Sigma}   
where \al{d_{IS,\Omega}(\Psi\|\Phi)=  \tr\int_{\Ts}\Omega(\log | \Psi | -\log |\Phi |+W_\Psi^{-1} \Phi W_\Psi^{-*}-I_m )} and $W_\Psi$ is a spectral factor of $\Psi$. It is not difficult to see that $d_{IS,\Omega}(\Psi\|\Phi)\geq 0$ and equality holds if and only if $\Phi=\Psi$. Moreover, $d_{IS,\Omega}(\Psi\|\Phi)$ is strictly convex with respect to $\Phi$. Notice that for $\Omega=I$ we obtain Itakura-Saito  distance. The Lagrangian of Problem (\ref{Pb_PSD_IS}) is 
\al{\Lc(\Phi,\Lambda)=\tr&\left[\int_{\Ts}\Omega(\log |\Psi | -\log |\Phi |+W_\Psi^{-1} \Phi W_\Psi^{-*}-I_m )\right.\nn\\ &\left.+G^*\Lambda G \Phi\right]-\tr(\hat \Sigma \Lambda)} which is bounded below only if $W_\Psi^{-*}\Omega W_\Psi^{-1}+G^*\Lambda G>0$. Moreover, it is not difficult to prove that $\Lc$ is strictly convex with respect to $\Phi$ and its point of minimum is given by setting equal to zero the first variation for any $\delta \Phi\in\Ls_\infty^{m\times m}(\Ts)$
\al{\delta \Lc(\Phi,\Lambda; \delta \Phi)=\tr&\left[\int_{\Ts}\Omega W_\Psi^{-1} \delta \Phi W_\Psi^{-*} +G^*\Lambda G \delta \Phi\right. \nn\\ &\left. -\Omega \int_0^\infty(\Phi+tI)^{-1}\delta\Phi(\Phi+tI)^{-1})\mathrm d t \right]\nn} leading to the optimality condition
 \al{\label{opt_IS}\int_0^\infty(\Phi+tI)^{-1}\Omega (\Phi+tI)^{-1}\mathrm d t= W_\Psi^{-*}\Omega W_\Psi^{-1}+G^*\Lambda G.} Clearly, it is not possible to find an explicit form for the optimal $\Phi$ and thus the dual analysis cannot carried out. Accordingly, we are not able to solve Problem (\ref{Pb_PSD_IS}). On the other hand, in the special case that $\Omega=\omega I$ with $\omega\in\Sc_1^+$, condition (\ref{opt_IS}) becomes 
 \al{\label{optPHi_IS_diag}\Phi=( \Psi^{-1}+\omega^{-1}G^*\Lambda G)^{-1}}
and the dual problem is equivalent to solve 
\al{\label{dual_IS}\Lambda^\circ= &\underset{\Lambda}{\mathrm{argmin\,}} -\int_\Ts \omega\log| \Psi^{-1}+\omega^{-1}G^*\Lambda G|+\tr(\hat\Sigma\Lambda)\nn\\  &\hbox{ s.t. } \Psi^{-1}+\omega^{-1}G^*\Lambda G>0.}
Moreover, it is possible to prove that Problem (\ref{dual_IS})  admits a solution and the solution to the primal problem is (\ref{optPHi_IS_diag}) with $\Lambda=\Lambda^\circ$.

\section{SIMULATION RESULTS}\label{sect_results}
We test the spectral estimator with indirect measurements using the weighted transportation distance. More precisely, we consider a Monte Carlo study constituted by 50 experiments. In each experiment:
\begin{itemize}
\item We generate a bivariate process $\xib$ with spectral density $\Phi_{\xib}=W_{\xib}W^*_{\xib}\in\Sc_2^+$, $W_{\xib}(z)=C_{\xib}(zI-A_{\xib})^{-1}B_{\xib}+D_{\xib}$ and $A_{\xib}$, $B_{\xib}$, $C_{\xib}$, $D_{\xib}$ are randomly generated. In particular, $A_{\xib}$ is a $4\times 4$ matrix whose eigenvalues are with absolute value less than or equal to 0.8. 
\item We generate the transfer matrix characteristic $H(z)^{-1}=C_H(zI-A_H)^{-1}B_H+D_H$. Matrices $A_H$, $B_H$, $C_H$, $D_H$ are randomly generated and such that $H(z)^{-1}$ is causally invertible. Moreover, $A_H$ is $4\times 4$ matrix whose eigenvalues are with absolute value less than or equal to 0.7. 
\item We generate $\Psi_{\xib}=W_pW_p^*$, $W_p(z)=C_p(zI-A_p)^{-1}B_p+D_p$, $A_p=A_{\xib}+\delta A$, $B_p=B_{\xib}+\delta B$, $C_p=C_{\xib}+\delta C$, $D_p=D_{\xib}+\delta D$. The perturbations matrices $\delta A$, $\delta B$, $\delta C$, $\delta D$ are randomly generated such that: its infinity matrix norm is equal to 0.08, $A_p$ is Schur stable and $\Psi_{\xib}\in \Sc_2^+$. 
\item We generate a short dataset $y(1)\ldots y(N)$, with $N=100$, extracted form a realization of $\yb$ whose shaping filter is $H^{-1}W_{\xib}$.  
\item We set the bank of filters as $G(z)=[\,z^{-l}I_m \, \ldots \,z^{-1}I_m\,]^T$ with $l=12$.
\item We compute $\hat \Sigma$ from the data using the optimization procedure in \cite{5953479} which guarantees the feasibility condition. 
\item We compute $\hat \Phi_{\xib,T}$ that is the solution to (\ref{Pb_PSD_orig}).
\item We compute the errors for each entry of the spectral density:
\al{\label{def_error}e_{\xib,T}^{(i,k)}(\vartheta)=| [\Phi_{\xib}(e^{j\vartheta})]_{i,k}-[\hat \Phi_{\xib,T}(e^{j\vartheta})]_{i,k}|.}
\end{itemize}  
In Figure \ref{fig_ris} we compare the average error defined in (\ref{def_error}) as a function of $\vartheta$ in the Monte Carlo study and the one obtained using the estimator $\hat \Phi_{\xib ,IS}$ defined in (\ref{Pb_PSD_origIS2}) which uses the Itakura-Saito distance. {\mz Furthermore, the average value of the $L_2$ norm of the error between the estimated spectral density and $\Phi_\xi$ is equal to $30.58$ using $\hat \Phi_{\xi,T}$ and $37.35$ using $\hat \Phi_{\xi, IS}$.} As we {\mz can} see, the estimator equipped with the transportation distance performs better than the one with the Itakura-Saito distance. Since the McMillan degree of $\hat \Phi_{\xib,T}$ is larger than the one of $\hat \Phi_{\xib,IS}$, as sanity check, we also considered $\hat \Phi_{\xib,IS}$ with larger values of $l$: we have obtained worse results than the ones with $l=12$. We conclude that the superiority of $\hat \Phi_{\xib,T}$ is probably due by the fact that the corresponding optimization scheme is adapted according to the  transfer matrix characteristic modeling the indirect measurements.

{\mz Finally, it is worth noting that the case of direct measurements, i.e. $H(z)=I$, has been analyzed in \cite{FERRANTE_TIME_AND_SPECTRAL_2012}. It has been shown that the estimator with the Itakura-Saito distance is better than the one with the Hellinger distance in the case that the shaping filter has some poles close to the unit circle.} 

\section{CONCLUSIONS}\label{sect_concl}
In this technical note, we have introduced the weighted transportation distance between Gaussian stationary processes. We showed that the latter corresponds to a weighted version of the Hellinger distance between multivariate power spectral densities. Then, we have developed a spectral estimator with indirect measurements using the transportation distance. Finally, a simulation study showed that the proposed estimator performs better than the one using the Itakura-Saito distance. {\mz It is worth noting that non-Gaussian processes as well as non-stationary processes are characterized by their joint probability density. Accordingly, the proposed optimal transport formulation seems the natural way for a possible extension in the aforementioned scenarios.}

\end{document}